\documentclass[onecolumn]{elsart3p}

\usepackage{amssymb,color}

\usepackage[english,francais]{babel}

\newtheorem{theorem}{Theorem}[section]
\newtheorem{lemma}[theorem]{Lemma}
\newtheorem{e-proposition}[theorem]{Proposition}
\newtheorem{corollary}[theorem]{Corollary}
\newtheorem{e-definition}[theorem]{Definition}
\newtheorem{remark}[theorem]{\it Remark\/}

\newtheorem{formal}[theorem]{Formal computation}

\newtheorem{theoreme}{Th\'eor\`eme}

\newtheorem{formel}{Calcul formel}

\newcommand{\R}{\mathbb{R}}
\newcommand{\Dt}{\partial_t}
\newcommand{\bea}{\begin{eqnarray}} 
\newcommand{\eea}{\end{eqnarray}}
\newcommand{\dx}{\,{\rm d}x}

\newcommand{\dz}{\,{\rm d}z}
\newcommand{\dt}{\,{\rm d}t}

\newcommand{\rf}[1]{(\ref{#1})}

\newcommand{\be}{\begin{equation}} 
\newcommand{\ee}{\end{equation}}
\newenvironment{multline*}{\begin{eqnarray*}}{\end{eqnarray*}}
\newenvironment{equation*}{\[}{\]}

\journal{the Acad\'emie des sciences}

\begin{document}

\centerline{Partial Differential Equations}
\begin{frontmatter}

\selectlanguage{english}

%%%%%%%%%%%%%%%%%%%%%%%%%%%%%%%%%%%%%%%%%%%%%%%%%%%%%%%%%%%%%%%%%%%%%%%%
%%%%%%%%%%%%%%%%%%%%%%%%%%%%%%%%%%%%%%%%%%%%%%%%%%%%%%%%%%%%%%%%%%%%%%%%

\title{Barenblatt profiles for a nonlocal porous medium equation}

\selectlanguage{english}

\author[Wroclaw]{Piotr Biler}, \ead{biler@math.uni.wroc.pl}
\author[Ceremade]{Cyril Imbert}, \ead{imbert@ceremade.dauphine.fr}
\author[Wroclaw]{Grzegorz Karch} \ead{karch@math.uni.wroc.pl}

\address[Wroclaw]{Instytut Matematyczny, Uniwersytet Wroc{\l}awski,
  pl. Grunwaldzki 2/4, 50--384 Wroc{\l}aw, Poland}
\address[Ceremade]{Universit\'e Paris-Dauphine, CEREMADE (UMR CNRS
  7534), Place de Lattre de Tassigny, 75775 Paris Cedex 16, France}

\medskip
\begin{center}
{\small Received ... \ ; accepted after revision ... \ .
\\ Presented by ...}
\end{center}

\begin{abstract}
  \selectlanguage{english} We study a generalization of the porous
  medium equation involving nonlocal terms.  More precisely, explicit
  self-similar solutions with compact support generalizing the
  Barenblatt solutions are constructed.  We also present a formal
  argument to get the $L^p$ decay of weak solutions of the corresponding Cauchy
problem. 

\vskip 0.5\baselineskip

\selectlanguage{francais}
\noindent{\bf R\'esum\'e} \vskip 0.5\baselineskip 
\noindent {\bf Solutions auto-similaires pour une \'equation des
  milieux poreux non locale}

Cette note est consacr\'ee \`a l'\'etude d'une g\'en\'eralisation non
locale de l'\'equation des milieux poreux. Plus pr\'ecis\'ement, on obtient
des formules explicites de solutions auto-similaires \`a support 
compact qui ressemblent fortement aux solutions de type
 Barenblatt.  On donne aussi un argument formel qui permet
 d'obtenir des estimations $L^p$ des solutions faibles du probl\`eme de Cauchy.
\end{abstract}

\end{frontmatter}

%%%%%%%%%%%%%%%%%%%%%%%%%%%%%%%%%%%%%%%%%%%%%%%%%%%%%%%%%%%%%%%%%%%%%%%%
%%%%%%%%%%%%%%%%%%%%%%%%%%%%%%%%%%%%%%%%%%%%%%%%%%%%%%%%%%%%%%%%%%%%%%%%
\selectlanguage{francais}
\section*{Version fran\c{c}aise abr\'eg\'ee}

Nous consid\'erons le probl\`eme de Cauchy pour l'\'equation non locale suivante
\be
\partial_tu-\nabla\cdot\left(u\nabla^{\alpha-1}\left(|u|^{m-1}\right)\right)=0,\label{eqm-f}
\ee 
avec $m> 1$, $x\in\R^d$, $\alpha\in(0,2)$, $t>0$, \`a laquelle on
ajoute une condition initiale  $u(0,x)=u_0(x)$. 
Ici $\nabla^{\beta}$ est un op\'erateur int\'egral
singulier g\'en\'eralisant le gradient usuel ($\beta=1$) et li\'e au
laplacien fractionnaire.

Le r\'esultat principal de cette note sont des formules explicites
de solutions auto-similaires qui se propagent \`a une vitesse
finie. 

%-----------------------------------------------------------------------
\medskip
%-----------------------------------------------------------------------
\begin{theoreme}[Solutions auto-similaires]\label{Thm:Main2F} \

\par\noindent 
La fonction 
$$
u(t,x)=Ct^{-\frac{d}{d(m-1)+\alpha}}\left(\left(R^2-|x|^2
    t^{-\frac{2}{d(m-1)+\alpha}}\right)_+^{\frac\alpha2}\right)^{\frac{1}{m-1}}
$$
(ainsi que ses translations en $x$) est une solution auto-similaire de
l'\'equation \rf{eqm-f} pour $R>0$ quelconque et une constante $C>0$
convenable.
\end{theoreme}
%-----------------------------------------------------------------------
%\bigskip
\medskip

Ensuite, nous pr\'esentons un calcul formel qui permet d'obtenir des
estimations en norme $L^p(\R^d)$ des solutions faibles du probl\`eme de
Cauchy, en particulier celles construites par Caffarelli et V\'azquez \cite{cv} dans le cas $m=2$
et telle que $|u_0(x)| \le C e^{-c|x|}$ pour deux
constantes $C$ et $c$.

\medskip

%-----------------------------------------------------------------------
\begin{formel}[Asymptotique du probl\`eme de Cauchy]\label{Thm:Main1F}\ 

\par\noindent
\'Etant donn\'ee une fonction $0\le u_0\in L^1(\R^d)\cap
L^\infty(\R^d)$, les normes $L^p(\R^d)$, $1\le p< \infty$, des
solutions faibles $u$  tendent
vers $0$ quand $t\to \infty$ avec le taux alg\'ebrique suivant
$$
\|u(t)\|_p\le C(d,\alpha,m,p)\|u_0\|_1^{\frac{\alpha + d(m-1)/p}{\alpha + d(m-1)}}
t^{-\frac{d}{d(m-1)+\alpha} (1-1/p)}\ \ \ {\rm for\ all}\ \ t>0.
$$ 
\end{formel}
%%%%%%%%%%%%%%%%%%%%%%%%%%%%%%%%%%%%%%%%%%%%%%%%%%%%%%%%%%%%%%%%%%%%%%

\selectlanguage{english}
\setcounter{equation}{0}
%%%%%%%%%%%%%%%%%%%%%%%%%%%%%%%%%%%%%%%%%%%%%%%%%%%%%%%%%%%%%%%%%%%%%%%%

\section{Introduction}\label{Sec:Intro}

We study a nonlocal generalization of the porous medium equation 
\be
\partial_tu-\nabla\cdot\left(u\nabla^{\alpha-1}\left(|u|^{m-1}\right)\right)=0,\label{eqm}
\ee 
where $m> 1$, $\alpha\in(0,2)$, $x\in\R^d$, $t>0$, supplemented with
an initial condition
\begin{equation}\label{eq:ic}
u(0,x) = u_0 (x). 
\end{equation}

The pseudodifferential (vector-valued) operator $\nabla^\beta$ in
\rf{eqm} is defined via the Fourier transform as $\nabla^\beta u=
{\mathcal F}^{-1}(i\xi|\xi|^{\beta-1} {\mathcal F}u)$.  This
definition is consistent with the usual gradient: $\nabla^1=\nabla$;
the components of $\nabla^0$ are the Riesz transforms; moreover we
have $\nabla\cdot\nabla^{\alpha-1}
=\nabla^{\frac\alpha2}\cdot\nabla^{\frac\alpha2}
=-(-\Delta)^{\frac\alpha2}$, where $(-\Delta)^{\frac\alpha2}$ denotes
the fractional Laplace operator:
$(-\Delta)^{\frac\alpha2}u={\mathcal F}^{-1}(|\xi|^\alpha{\mathcal
  F}u)$.  It can also be defined by real analysis tools as follows
$\nabla^{\alpha-1} u = \nabla I_{2-\alpha}u$, where $I_\beta$ for
$\beta\in(0,d)$ is an integral smoothing operator, called the Riesz
potential (see \cite[Ch. V]{S})
$$
I_\beta(u)(x)=-C_\beta\int \frac{u(x+z)}{|z|^{d-\beta}}\dz 
$$
with some $C_\beta>0$. Note that then $\nabla^{\alpha-1} u(x)=\nabla
I_{2-\alpha}(u)(x)=(d+\alpha-2) C_{2-\alpha}\int
(u(x+z)-u(x))\frac{z}{|z|^{d+\alpha}}\,{\rm d}z$, $\alpha\in(0,2)$.

Eq.~\rf{eqm} can be interpreted as a transport equation of the type
$\partial_t u = \nabla \cdot ( u \mathbf{v})$ for some velocity vectorfield
$\mathbf{v}$ which is a potential; more precisely, $\mathbf{v}=\nabla
\mathbf{p}$ where ${\mathbf p}=I_{2-\alpha}(|u|^{m-1})$. It can be interpreted as
a nonlocal pressure in the case of nonnegative initial data. Then, of course
$(-\Delta)^{\frac{2-\alpha}{2}}{\mathbf p}=|u|^{m-1}$, see \cite{cv},
\cite{cv2} for that notation.

Notice that for $\alpha=2$ and nonnegative initial data we recover the
Boussinesq equation ($m=2$), and the usual porous media equation
($m>1$): $\Dt u = \nabla \cdot \left(u \nabla
  \left(u^{m-1}\right)\right )$, $t>0$, $x\in \R^d.$

Recently, L. Caffarelli and J. L. V\'azquez (\cite{cv}, \cite{cv2}) studied
equation \rf{eqm} in the case $m=2$. They proved the existence of weak
solutions for nonnegative bounded integrable initial data with exponential
decay at infinity.  They also treat the case of bounded and compactly
supported initial data, which propagate with finite speed. It is shown in
\cite{cv} that self-similar solutions can be constructed by considering an
obstacle problem for the fractional Laplace operator.  In this note, we
contribute to those results constructing {\em explicit compactly supported
  self-similar solutions}.  Moreover, we show a kind of hypercontractivity
estimates, \textit{i.e.} the optimal decay in $L^p$ 
of general solutions with $u_0\in L^1(\R^d)$.

Results on equation~\rf{eqm} in this note are multidimensional
generalizations of those obtained in \cite{bkm} for a~model of the
dynamics of dislocations in crystals (for the integral of $u$ when
$d=1$ and $m=2$).  The structure of \rf{eqm} suggests that it should
enjoy the conservation of mass and some comparison properties as was
shown in \cite{bkm}.
For an analysis of a  related nonlocal equation, see \cite{im}. 

Complete proofs of all results announced in this note will be published in \cite{fpm}.

\section{Self-similar solutions}
\label{sec:self-sim} 

The equation \rf{eqm} has the following scaling property: 
$$
{\rm if\ } u(t, x)\ {\rm is\ a\ solution,\ so\ is\ }
\ell^{d\lambda}u(\ell t,\ell^\lambda x)
$$ 
for each $\ell>0$ and $\lambda=\frac{1}{d(m-1)+\alpha}$.  We look for
nonnegative solutions that are invariant under that scaling, \textit{i.e.}  of
the following form 
\be
\label{eq:self-sim-shape} 
u(t,x) =
\frac{1}{t^{d\lambda}}\Phi_{\alpha,m}\left(\frac{x}{t^{\lambda}}\right),\
\ {\rm where}\ \ \ \lambda=\frac{1}{d(m-1)+\alpha}, 
\ee 
for a function $\Phi_{\alpha,m}: \R^d\to \R^+$ satisfying the
following nonlinear and nonlocal equation in $\R^d$ 
\be
\label{eq:ellip-sim} 
-\lambda \nabla \cdot (y \Phi_{\alpha,m}) = \nabla \cdot
\left(\Phi_{\alpha,m} \nabla^{\alpha-1} \Phi_{\alpha,m}^{m-1}\right).
\ee 
Before stating our main result, we recall the definition of weak
solutions for \rf{eqm} introduced in \cite{cv} in the case $m=2$.
\medskip

\begin{e-definition}[Weak solutions]
A function $u: (0,T) \times \R^d \to \R$
is a \emph{weak solution} of \rf{eqm} in $Q_T= (0,T) \times \R^d$
submitted to the initial condition $u(0,x) = u_0(x)$ if $u \in L^1
(Q_T)$, $I_{2-\alpha} (|u|^{m-1}) \in L^1 (0,T;W^{1,1}_{\rm {loc}}
(\R^d))$, $u \nabla I_{2-\alpha} (|u|^{m-1}) \in L^1 (Q_T)$ and \be
\int\int u (\varphi_t - \nabla I_{2-\alpha} (|u|^{m-1}) \cdot \nabla
\varphi) \,{\rm d}x \,{\rm d}t + \int u_0 (x) \varphi (0,x) \, {\rm
  d}x = 0
\label{weak}
\ee for each test function $\varphi \in C^1 (Q_T)$ such that
 $\varphi$ has compact support in the
space variable $x$, and vanishes near $t=T$.
\end{e-definition}
\medskip

%----------------------------------------------------------------
\begin{theorem}[Self-similar solutions]\label{thm:self-sim} 
  For each $\alpha\in(0,2]$, $m>1$ and $R>0$, the function \be
  \Phi_{\alpha,m} (y) =
  \left(k(R^2-|y|^2)_+^{\frac{\alpha}{2}}\right)^{\frac{1}{m-1}}\ \
  {\rm with}\ \ \
  k=\left(\frac{d}{d(m-1)+\alpha}\right)\left(\frac{\Gamma\left(\frac{d}{2}\right)}{2^\alpha\Gamma\left(1+\frac{\alpha}{2}\right)\Gamma\left(\frac{d+\alpha}{2}\right)}\right)
\label{Phi}
\ee
is a solution of \rf{eq:ellip-sim}. Consequently, the function 
\be 
u(t,x)=t^{-\frac{d}{d(m-1)+\alpha}} \left(k 
  \left(R^2-|x|^2t^{-\frac{2}{d(m-1)+\alpha}}\right)_
  +^{\frac\alpha2}\right)^{\frac{1}{m-1}}
\label{self}
\ee
is a weak solution of \rf{eqm}, satisfies the equation in the
pointwise sense for $|x|\neq Rt^{\frac{1}{d(m-1)+\alpha}}$, and is
\newline 
$\min\left\{\frac{\alpha}{2(m-1)},1\right\}$-H\"older continuous at the interface
$|x|=R t^{\frac{1}{d(m-1)+\alpha}}$.
\end{theorem}
%-----------------------------------------------------------------
\medskip

When $\alpha=2$, we recover the classical Kompaneets--Zel'dovich--{\bf
  Barenblatt}--Pattle formulas, see, \textit{e.g.}, \cite{V2}.

Note also that given $M>0$ there exists a unique $R>0$ such that
$\int\Phi_{\alpha,m}(y)\,{\rm d}y=M$.  \medskip

\begin{remark}
As mentioned above, self-similar solutions of \rf{eqm} have been
proved to exist in \cite{cv} by studying the following obstacle
problem for the fractional Laplacian
 $$
 P\ge \Phi,\ \ \ V=(-\Delta)^{\frac{\alpha}{2}}P\ge 0,\ \ 
 {\rm either\ \ } P=\Phi\ \  {\rm or\ \ } V=0,
$$ 
with $\alpha\in(0,2)$ and $\Phi(y)=C-a|y|^2$.  The novelty of our
approach is that we exhibit the explicit solution of this obstacle
problem: $P(y)=I_\alpha\left(\Phi_{\alpha,2}\right)(y)$, where
$I_\alpha$ is the Riesz potential, and $\Phi_{\alpha,2}$ is defined in
\rf{Phi} with $m=2$ and a suitable $R>0$.
\end{remark}

Those explicit self-similar solutions express one of the most
important features of the porous medium equation: the property of
finite propagation speed.  In the case of the classical porous medium
equation ($\alpha=2$), this property is established using
comparison with suitably large self-similar solutions,
\textit{cf.} \cite{V2}. For the generalized porous medium equation \rf{eqm}
with $m=2$, special supersolutions have been used for comparison,
\textit{cf.} \cite{cv2}.
\medskip

The proof of Theorem~\ref{thm:self-sim} is based on an application of
the following fundamental technical fact.  \medskip

%------------------------------------------------------
\begin{lemma}\label{lem:getoor}
For all $\beta\in(0,2)$ and $\gamma>0$, we have 
\begin{equation}\label{eq:getoor}
I_\beta\left((1-|y|^2)_+^{\frac{\gamma}{2}}\right)= 
\left\{
  \begin{array}{ll} 
    C_{\gamma,\beta,d} \times
   {}_2F_1\left(\frac{d-\beta}{2},-\frac{\gamma+\beta}{2}; \frac{d}{2};|y|^2\right)\ & \quad  {\rm for\ \ } |y| \le 1, \\
    \tilde{C}_{\gamma,\beta,d} |y|^{\beta-d} 
   \times {}_2F_1\left(\frac{d-\beta}{2},\frac{2-\beta}{2}; \frac{d+\gamma}{2};\frac{1}{|y|^2} \right) 
    & \ \ \ {\rm for\ \  } |y| \ge 1,
  \end{array}
\right.
\end{equation}
with $C_{\gamma,\beta,d} = 2^{-\beta} \frac{\Gamma
  \left(\frac\gamma2+1\right)\Gamma
  \left(\frac{d-\beta}2\right)}{\Gamma \left(\frac{d}2\right)\Gamma
  \left(\frac{\beta+\gamma}2+1\right)}$ and  
$\tilde{C}_{\gamma,\beta,d} = 2^{-\beta} \frac{\Gamma
  \left(\frac\gamma2+1\right)\Gamma
  \left(\frac{d-\beta}2\right)}{\Gamma \left(\frac{d}4\right)\Gamma
  \left(\frac{d+\gamma}2+1\right)}$,  where ${}_2F_1$ is the hypergeometric function.
\end{lemma}
\medskip

The verification of \rf{eq:getoor} consists in passing to Fourier
transforms and calculating certain integrals (the so called
(Sonine--)Weber--Schafheitlin integrals) involving Bessel and
hypergeometric functions ${}_2F_1$, \textit{cf.} \cite{mos}.  \medskip

\noindent{\em Proof of Theorem~\ref{thm:self-sim}.} \ \ \ Let
$\phi_\alpha(y)=\left(1-|y|^2\right)_+^{\frac{\alpha}{2}}$.  Observe
that $\phi_\alpha\in L^1(\R^d)$. We next show that $I_{2-\alpha} (\phi_\alpha) \in
W^{1,1}_{\rm {loc}}(\R^d)$.  By Lemma~\ref{lem:getoor} with
$\gamma = \alpha \in (0,2)$, $\beta = 2 -\alpha$, and
${}_2F_1(a,-1;c;z)= 1 - \frac{a}c z$, we get
$$
I_{2-\alpha} (\phi_\alpha) (y) = \left\{\begin{array}{ll}
    C_{\alpha,2-\alpha,d} \left( 1 - \frac{d+\alpha-2}{d} |y|^2 \right) & \ \ {\rm if \ \ } |y| \le 1,\\
    \tilde{C}_{\alpha,2-\alpha,d} |y|^{2-(d+\alpha)} {}_2F_1
    \left(\frac{d+\alpha}2 -1,\frac\alpha2;\frac{d+\alpha}2;
    \frac{1}{|y|^2}\right) & \ \ {\rm if \ \ } |y| \ge 1,
\end{array}\right.
$$
which is a locally integrable function. Recalling that
$\nabla^{\alpha-1} = \nabla I_{2-\alpha}$, we then deduce by the chain
rule that for $y \in B_1$, $\nabla^{\alpha -1} (\phi_\alpha) (y)
=-\left(d\,K_{\alpha,d}\right)^{-1} y$ where $K_{\alpha,d}$ is defined in Corollary~\ref{lem:getoor-genuine} 
below.
For $|y|\ge 1$, one uses $\frac{\partial}{\partial z}\,{}_2 F_1
(a,b;c;z)= \frac{ab}{c} \; {}_2 F_1 (a+1,b+1;c+1;z)$, hence
$\nabla^{\alpha-1} (\phi_\alpha)$ is locally integrable.  To conclude,
we remark that
$$
\phi_\alpha^{\frac{1}{m-1}} (y) \nabla^{\alpha-1} \phi_\alpha (y)=
-\phi_\alpha^{\frac{1}{m-1}}(y) \,{(K_{\alpha,d} d)^{-1}}\, y, \ \
{\rm for\ all}\ \ y \in \R^d,
$$
because $\phi_\alpha(y)=0$ for $|y|\ge 1$. 

Hence,  scaling the variables, we immediately obtain that the function $\Phi_{\alpha,m}$ defined in \rf{Phi}
satisfies  $\nabla^{\alpha-1}(\Phi_{\alpha,m}^{m-1})=-\lambda y$ for $|y|< R$. 
Now, for all $y\in\R^d$ the identity 
$$ 
-\lambda y \Phi_{\alpha,m} =   \Phi_{\alpha,m} \nabla^{\alpha-1} \Phi^{m-1}_{\alpha,m}
$$
follows because $\Phi_{\alpha,m}(y)=0$ for all $|y|\ge R$, so
\rf{eq:ellip-sim} holds with $\lambda=\frac{1}{d(m-1)+\alpha}$ and
$k=d\lambda\, K_{\alpha,d}$ defined in \rf{Phi}.

It is straightforward to verify using \rf{eq:ellip-sim} that $u$ given
by formula \rf{self} is a weak solution of \rf{eqm} in each strip
$(t_0,T)\times\R^d$, $0<t_0<T<\infty$. Moreover, the H\"older continuity is easy to check.
\hfill$\square$
 
The following known result (with an important probabilistic
interpretation) proved by Getoor \cite[Th. 5.2]{G}, see also
\cite[App.]{L} for a related calculation, is an immediate consequence
of Lemma \ref{lem:getoor}.  For its proof, it suffices to use the
relation $(-\Delta)^{\frac\alpha2} = \Delta I_{2-\alpha}$.
\medskip

\begin{corollary}\label{lem:getoor-genuine} {For each $\alpha \in
    (0,2]$, the identity $K_{\alpha,d}(-\Delta)^{\frac\alpha2}\left(
      \left(1-|y|^2\right)_+^{\frac{\alpha}{2}} \right)=
    -1\ \ {\rm in}\ \ B_1 $ holds with the explicit
    constant $K_{\alpha,d}=\frac{\Gamma\left(\frac{d}2\right)
    }{2^\alpha \Gamma\left(1+\frac\alpha2\right) \Gamma
      \left(\frac{d+\alpha}2\right)}$. }
\end{corollary} 
%-------------------------------------
\section{The Cauchy problem and asymptotics}\label{Sec:CP}

We now briefly discuss the questions of the existence of weak solutions and
their uniqueness in the case $m=2$.  In \cite{bkm}, viscosity solutions
have been considered which permitted the authors to prove regularity and
uniqueness of solutions to \rf{eqm} in one space dimension.  In higher
dimensions, a construction of mild solutions is achieved in \cite{cv}
through a parabolic regularization of \rf{eqm} with a suitable cutoff of
the singular kernel of $I_{2-\alpha}$, and then a passage to the limit. The
uniqueness of weak solutions and, {\em a fortiori}, the validity of the
full comparison principle seem to be difficult questions, {\it cf.} a
discussion in \cite{cv}. Another construction of weak solutions to \rf{eqm}
with $m>1$ as limits of mild solutions of parabolically perturbed equation
\rf{eqm} will be published in \cite{fpm}.  \medskip

%-------------------------------------------------------------------------
\begin{formal}[Decay of solutions for the Cauchy problem]
\label{thm:existence} 
\par\noindent
Suppose that for $0\le u_0\in L^1(\R^d)\cap L^\infty(\R^d)$.  Then
$L^p(\R^d)$ norms ($1\le p< \infty$) of any sufficiently regular and global
in time nonnegative weak solution $u$ of \rf{eqm}--\rf{eq:ic} such that
$\int u(t,x)\dx=\int u_0(x)\dx$ decay algebraically
\begin{equation}\label{decay}
\|u(t)\|_p\le C(d,\alpha,m,p)\|u_0\|_1^{\frac{\alpha + d(m-1)/p}{\alpha + d(m-1)}}
t^{-\frac{d}{d(m-1)+\alpha} (1-1/p)}\ \ \ {\rm for\ all}\ \ t>0.
\end{equation}
\end{formal}

%--------------------------------------------------------------------------

\medskip

Our computation below can be applied, for example, to weak solutions constructed  by 
Caffarelli and  V\'azquez \cite{cv}, who studied problem
\rf{eqm}--\rf{eq:ic} with $m=2$ and with initial conditions satisfying 
$
0 \le
u_0(x) \le C e^{-c|x|}$ for some constants $c,C$ and all $x\in \R^d$.

\medskip\noindent
{\it Formal proof of \rf{decay}.}
In order to prove the announced $L^p$ estimates of solutions (similar
to those for degenerated partial differential equations like the
porous medium equation in \textit{e.g.} \cite{V1}, \cite[Ch. 2]{CJMTU}), we
recall  the
Stroock--Varopoulos inequality for $q \ge 1$
\be \int |w|^{q-2}w(-\Delta)^{\frac\alpha2}w\dx\ge
\frac{4(q-1)}{q^2}\int\left|\nabla^{\frac\alpha2}|w|^{\frac{q}2}\right|^2\dx
\label{SV}
\ee 
valid for each $w\in L^q(\R^d)$ such that
$(-\Delta)^{\frac\alpha2}w\in L^q(\R^d)$.  Note that the constant in
\rf{SV} is the same as for the usual Laplacian operator $-\Delta$
(\textit{i.e.} $\alpha=2$). The proof is given, \textit{e.g.}, in \cite[Prop. 1.6]{LS}
and \cite[{Th. 2.1 combined with (1.7)}]{LS}.  

We will also need the following Nash inequality
\be 
\|v\|_2^{2(1+\frac\alpha{d})}\le
C_N\|\nabla^{\frac\alpha2}v\|^2_2\|v\|_1^{\frac{2\alpha}d}\label{N}
\ee 
valid for all functions $v$ with $v\in L^1(\R^d)$, 
$\nabla^{\frac\alpha2}v\in L^2(\R^d)$ 
with a constant $C_N=C(d,\alpha)$.  
The proof of \rf{N} for $d=1$ can be found in, \textit{e.g.}, \cite[Lemma 2.2]{KMX}, 
and this extends easily to the general case $d\ge 1$. 

Moreover, we will need the Gagliardo--Nirenberg type inequality: for $p>1$ and $p\ge m-1$, 
\be
\|u\|_p^a\le C_N\left\|\nabla^{\frac\alpha2}|u|^{\frac{r}2}\right\|^2_2\|u\|_1^b 
\label{inter}
\ee 
with $a=\frac{p}{p-1}\frac{d(r-1)+\alpha}{d}$ and 
$b=\frac{p\alpha +d(m-1)}{d(p-1)}$ and  $r=p+m-1$. 
This inequality is a consequence of the Nash inequality \rf{N} written 
for $v=|u|^{\frac{r}2}$, \textit{i.e.} 
$\|u\|_r^{r(1+\frac\alpha{d})}\le
C_N\left\|\nabla^{\frac\alpha2}|u|^{\frac{r}2}\right\|^2_2\|u\|_{\frac{r}2}^{\frac{r\alpha}{d}}$,  
and two H\"older inequalities for the $L^q$ norms: 
$ 
\|u\|_p\le \|u\|_r^\gamma\|u\|_1^{1-\gamma}$ and 
$\|u\|_{\frac{r}2}\le \|u\|_p^\delta\|u\|_1^{1-\delta},
$ 
which hold with $\gamma=\frac{p-1}{r-1}\frac{r}{p}$ and
$\delta=\frac{r-2}{p-1}\frac{p}{r}$. Combining the above three
inequalities, we get \rf{inter}. 
\medskip

%---------------------------------------------------------------
 The
computation, which lead to \rf{decay},  consists in getting a differential inequality of the form
$ \frac{\rm d}{\dt}\int|u|^p\dx\le - K\|u\|_p^a\|u\|_1^{-b} $ for some
positive constant $K$ and where $a$ and $b$ appear in (\ref{inter}).

Multiply \rf{eqm} by $u^{p-1}$ with $p>1$, integrate by parts, and use
the relation
$\nabla\cdot\nabla^{\alpha-1}=-(-\Delta)^{\frac{\alpha}{2}}$ to get
\bea \frac1p\frac{\rm d}{\dt}\int u^p\dx &=& -(p-1)\int u^{p-1}\nabla
u\cdot\nabla^{\alpha-1}(u^{m-1})\dx
\nonumber\\
&=& -\frac{p-1}{p}\int u^{p}(-\Delta)^{\frac\alpha2} u^{m-1}\dx
\nonumber\\&\le&
-\frac{4(p-1)(m-1)}{(p+m-1)^2}\left\|\nabla^{\frac\alpha2}\left(u^{\frac{p+m-1}2}\right)
\right\|^2_2\nonumber \eea after applying the Stroock--Varopoulos
inequality \rf{SV} with $w=u^{m-1}$ and $q=\frac{p}{m-1}+1$.  To
estimate the right hand side of the above inequality, we use
 \rf{inter}. We finally have:
$$ 
\frac{\rm d}{\dt}\int|u|^p\dx\le - K\|u\|_p^a\|u\|_1^{-b}
$$ 
with 
$K=\frac{1}{C_N}\frac{4p(p-1)(m-1)}{(p+m-1)^2}$. The above inequality
leads to the differential inequality $ \frac{\rm d}{\dt}f(t)\le
-KM^{-b}f(t)^{\frac{a}p} $ for the function $f(t)=\|u(t)\|_p^p$,
$M=\|u_0\|_1$, and $a/p>1$, which immediately gives the algebraic
decay of the $L^p$ norms for $p\ge m-1$: $f(t)\le
\left(K\left(\frac{a}p-1\right)M^{-b}\,
  t\right)^{-\frac{1}{\frac{a}p-1}}$. Finally, we obtain the desired
estimate with
$$
C(d,\alpha,m,p) = \left(K\left(\frac{a}{p}-1 \right)\right)^{-\frac1{a-p}}= \left[ \frac{4 (m-1)(d(m-1)+\alpha)}{C_N d} \times
  \frac{p}{(p+m-1)^2} \right]^{-\frac{d}{d(m-1)+\alpha} (1-\frac{1}p)} .
$$

 Of course, this is sufficient
to get the conclusion of Theorem \ref{thm:existence} for each $1\le p
<\infty$ since $\frac{1}{\frac{a}p-1}=\frac{d(p-1)}{d(m-1)+\alpha}$,
and interpolating between $L^1(\R^d)$ and $L^p(\R^d)$ with $p$
sufficiently large.  \hfill$\square$

\medskip

\begin{remark}In our recent work  \cite{fpm}, we present a
 more subtle iterative argument, which allows us to show that also 
$$\|u(t)\|_\infty\le C\|u_0\|_1^{\frac{\alpha}{d(m-1)+\alpha}}t^{-\frac{d}{d(m-1)+\alpha}}
\quad\mbox{for all}\quad  t>0.$$
\end{remark}

\medskip
\noindent {\bf Acknowledgments.\ } The preparation of this paper was
supported by the Polish Ministry of Science (\hbox{MNSzW}) grant N201
418839, and a PHC POLONIUM project 0185, 2009--2010, and the
Foundation for Polish Science operated within the Innovative Economy
Operational Programme 2007--2013 funded by European Regional
Development Fund (Ph.D. Programme: Mathematical Methods in Natural
Sciences).  The authors wish to thank Juan Luis V\'azquez, R\'egis
Monneau and Jean Dolbeault for fruitful discussions they had together.

%%%%%%%%%%%%%%%%%%%%%%%%%%%%%%%%%%%%%%%%%%%%%%%%%%%%%%%%%%%%%%%%%%%%%%%%

\end{document}